\documentclass[12pt,a4paper]{article}
\usepackage[greek, english]{babel}
\usepackage[utf8x]{inputenc}
\usepackage{amsmath}
\usepackage{amsfonts}
\usepackage{amssymb}
\usepackage{theorem}
\usepackage{makeidx}
\usepackage[pdftex,final]{graphicx}
\usepackage[dvipsnames,usenames]{color}
\DeclareGraphicsExtensions{.pdf, .jpg}
\usepackage{rotating}
\usepackage{float}
\usepackage{boxedminipage}
\usepackage{cancel}
\usepackage{a4wide}
\theorembodyfont{\upshape}
\bibliography{References}

\newcommand{\indt}{\indent}
\newcommand{\en}{\selectlanguage{english}}

\newcommand{\bd}{\textbf}
\newcommand{\mbd}{\mathbf}
\newcommand{\R}{\mathbb{R}}
\newcommand{\s}{\sigma}

\newtheorem{definition}{Definition}

\begin{document}
\en
\title{Metrics and norms used for obtaining sparse solutions to underdetermined Systems of Linear Equations}
\author{Leoni Dalla and George K. Papageorgiou} 
\maketitle
\begin{abstract}
This paper focuses on defining a measure, appropriate for obtaining optimally sparse solutions to underdetermined systems of linear equations.\footnote{The following work done, was within the completion of my master thesis titled ``Algorithms for the computation of sparse solutions of undefined systems of equations" at the department of Mathematics, University of Athens which was assigned to me in association with the department of Informatics and Telecommunications, National and Kapodistrian University of Athens.} The general idea is the extension of metrics in n-dimensional spaces via the Cartesian product of metric spaces. 
\end{abstract}

\section{Introduction}
In general topology, mathematicians have long ago defined measures that had seen minimum usage (if not at all) in applications. Later on, the development of Measure Theory was mandatory for the progression of applied mathematics and other sciences too. Along with the progress in computer sciences came the demand defining measures of unusual nature and uncovering the properties they obey.

In signal (image or sound) processing a usual problem that arises, is how to transfer a signal using a sparse (economical, but sufficient) representation \cite{Siam, donoho}. Given a specific matrix $\mbd{A}$ of dimension $m\times n,\ m,n\in \mathbb{N}$ with $m<n$ (underdetermined) and a vector $\mbd{b},$ find among all, the  sparsest or (a less sparse) solution $\mbd{x}\in \R^n$ of the linear system $\mbd{Ax=b}$. This is the simplest form of the problem, which means that the noise of the signal is not included (noiseless problem).

Since an undefined system of linear equations has infinite number of solutions, they need to be filtered, using additional functions, in order to obtain solutions of a certain type according to specific criteria. Functions that measure ``energy", like the $l_2$ norm, are used in many occasions, yet measuring sparsity requires a measure of ``sparsity", i.e. a different function \cite{Siam, donoho}. 

The optimization task is minimizing the number of nonzero coordinates of a vector in $n$-dimensions, i.e. finding a sparse representative of the signal. The number of nonzero coordinates of a vector $\mbd{x}$ is known to be the number of elements included in the set of nonzero values of a vector, which is called \text{\em support} of the vector, i.e. $supp\{\mbd{x}\}$. Also, in recent work of \text{\em Donoho} and \text{\em Elad} the measure was ``used" under the symbol of norm $\|\mbd{x}\|_0= \# \{ i: x_i\neq 0\},$ but it is clear that it does not satisfy the norm properties \cite{Siam, donoho}. 

In the following paper, we begin posing some examples in everyday life, where different measures are needed in order to figure out distances. After a short review on metric spaces follows the definition of $p$-metrics in a Cartesian product space. The next section is of main interest, since we equip $\R^n$ with metrics and prove that the discrete metric could be obtained as a limit of a $p$-metric \cite{Siam}. In addition, follows a review on norms and the correlation between norms and metrics. Finally, we conclude with a comparison of functions, on the quest for a convex one, suitable for optimization tasks.

\section{Measures in everyday life}
In everyday life people subconsciously use measures in order to figure out distances, albeit those measures are not always well defined. Given an arbitrary set of points, a matter of most concern is to measure the distances between those points. However, the measure that we use in every problem is different and depends on the scale we would like to use, as well as the structure of the setting.

\begin{figure}[htb]
\centering
\includegraphics[width=110mm, height=70mm]{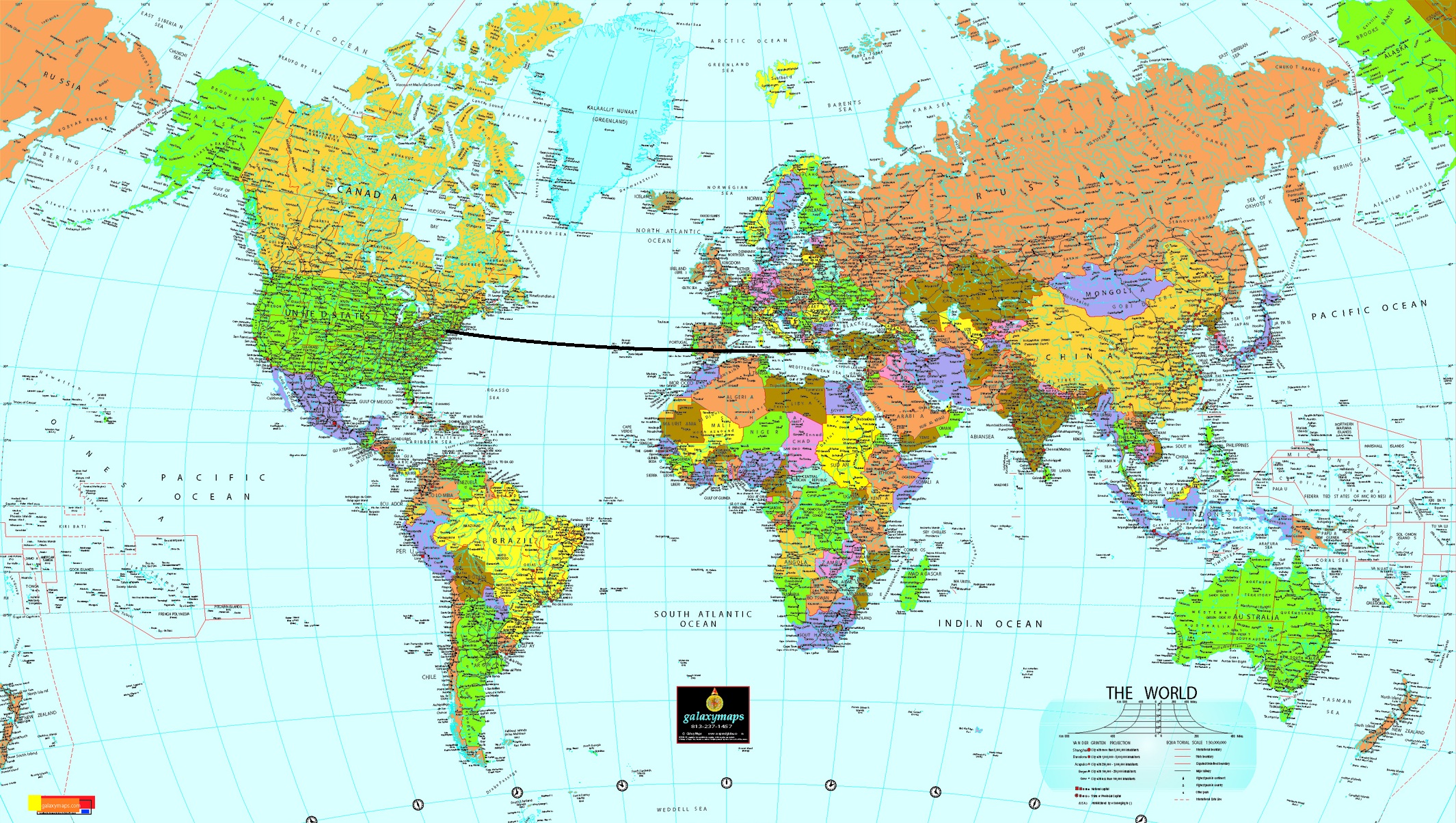}
\caption{\small{Geodesic distance between Athens - New York: $7920\ km$}}
\label{NY-A}
\end{figure}

The distance between Athens and New York is measured using the geodesic line between those points, i.e. the shortest route between two points on Earth's surface. 
(Fig.~\ref{NY-A}).
\begin{figure}[!]
\centering
\includegraphics[scale=0.17]{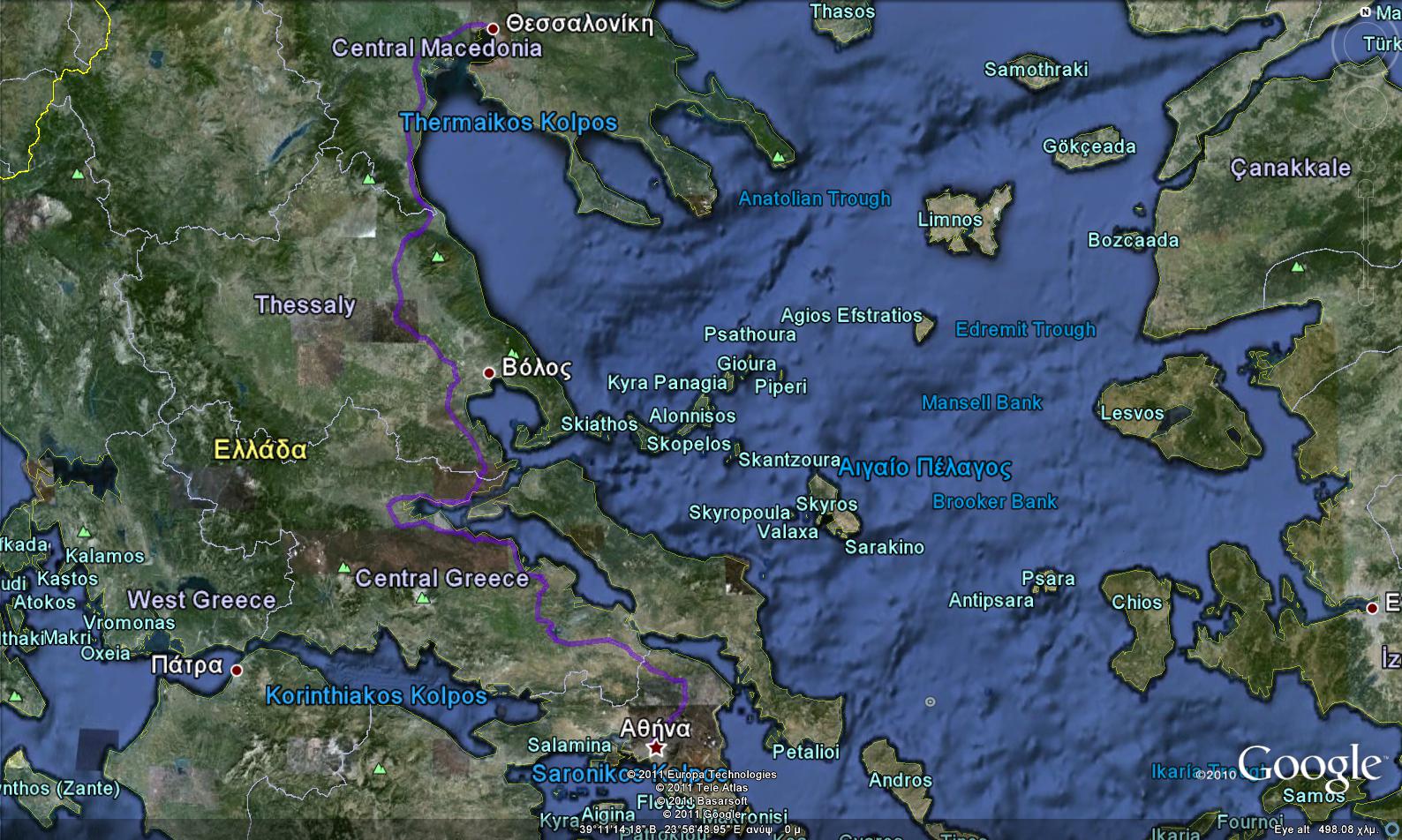}
\caption{\small{Travel distance between Athens - Thessaloniki: $502\ km$}}
\label{A-Th}
\end{figure}

\begin{figure}[!]
\centering
\includegraphics[ width=60mm, height= 85mm]{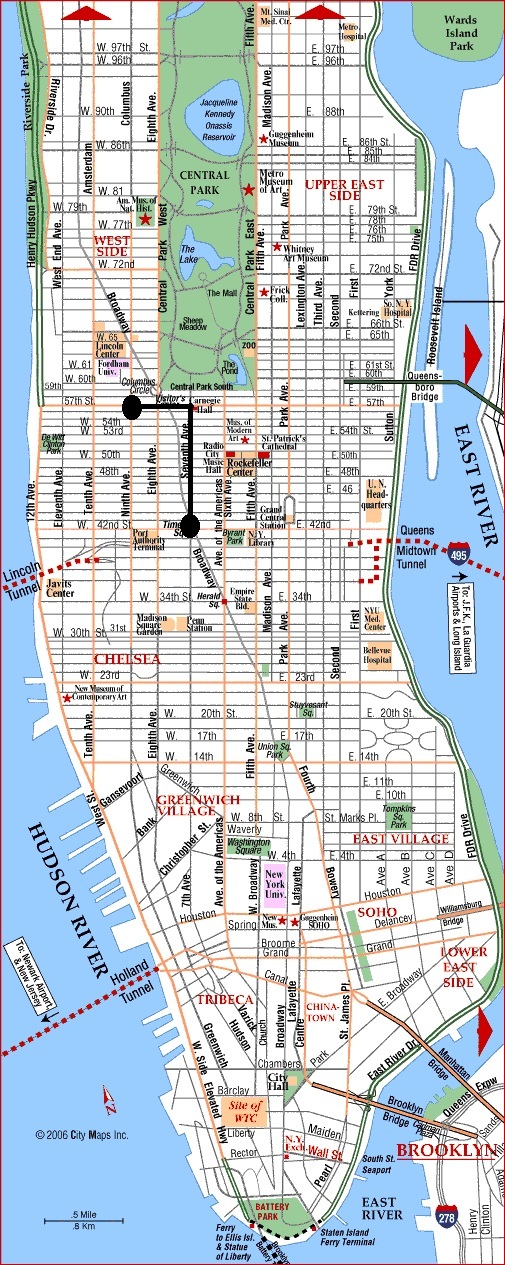}
\caption{\small{Distance in the area of Manhattan, borough of New York City.}}
\label{Man}
\end{figure} 
In case of a road trip, the travel distance depends on the road's structure and does not coincide with shortest distance between those two points (towns) (Fig.~\ref{A-Th}).

Furthermore, the existence of distances that differ from our perception of the shortest path cannot pass unnoticed. The distance that a person has to travel in the area of Manhattan (borough of New York City) in order to move from Times Square to the junction of 57th Street with 9th Avenue depends on the structure of the setting (Fig.~\ref{Man}).

\begin{figure}[htb]
\centering
\includegraphics[ scale=0.4]{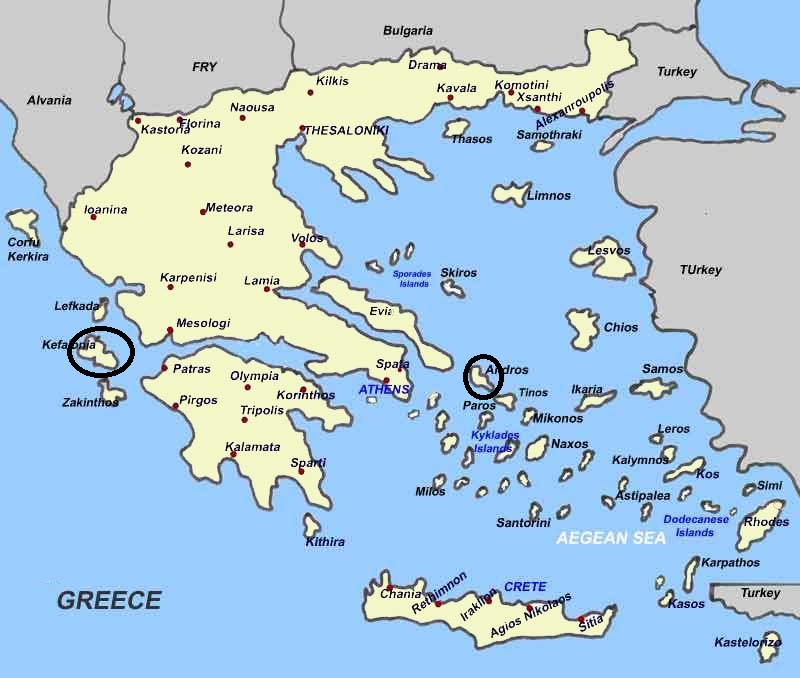}
\caption{\small{Distance between Kefallonia - Andros}}
\label{K-A}
\end{figure}

Another measure, used in order to define distances between compact sets, e.g. the distance between two islands, is the Hausdorff distance (named after Felix Hausdorff) between the whole sets $K$ and $A$ defined 
\[h(K,A)=max\{\max_{\kappa \in K}{ \min_{\alpha \in A}|\kappa - \alpha|}, \max_{\alpha \in A} \min_{\kappa \in K}|\kappa - \alpha| \}. \]
The latter represents, e.g. the minimum distance one has to travel in order to move from any village of the island of Andros (or Kefallonia) to any village of the island of Kefallonia (or Andros) (Fig.~\ref{K-A}). 

So far, we have seen cases where the concept of distance needs to be mathematically defined in order to understand, develop and solve problems arising from very different settings. Hence, we should define the means needed in order to measure in a wide variety of cases. 

\section{Metric spaces}

\begin{figure}[htb]
\centering%
\includegraphics[scale=0.4]{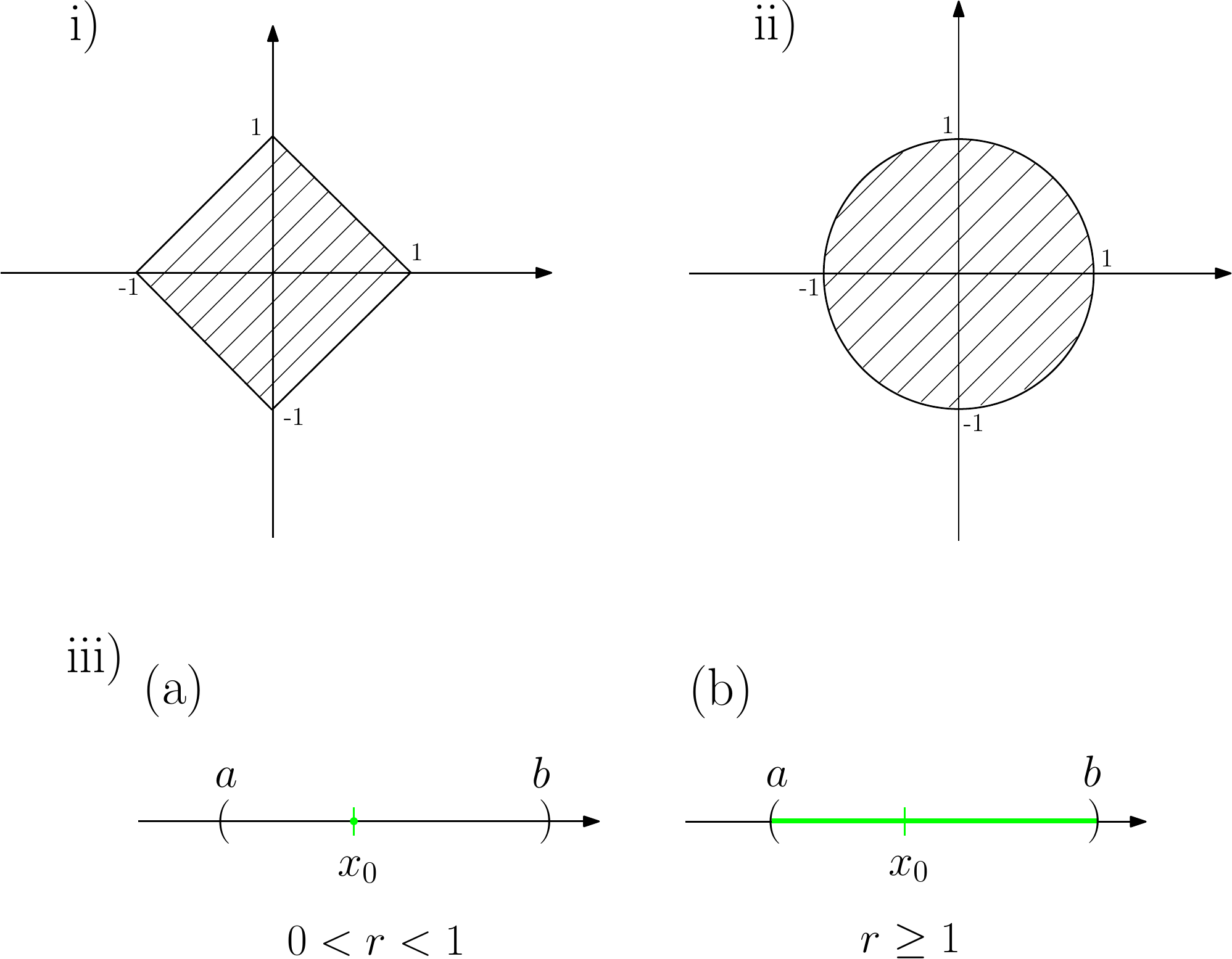}
\caption{\small{i) $ \widetilde{S}_{d_1}(\mbd{0},1)=\{\mbd{x}\in \R^2: d_1(\mbd{x},\mbd{0})\leq 1 \}.$ ii) The closed ball $\widetilde{S}_{d_2}(\mbd{0},1)=\{\mbd{x}\in \R^2: d_2(\mbd{x},\mbd{0})\leq 1 \}$ of the Euclidean metric in $\R^2$ is the unit circle. iii) ($a$) For $x_0\in X=(a,b)\subseteq \R$ and $0<r<1,$  $\widetilde{S}_{\s_0}(x_0, r)=\{x_0 \}.$  iii) ($b$) For $x_0\in X=(a,b)\subseteq \R$ and $r\geq 1$ the closed ball is the entire set $X,$ i.e.  $\widetilde{S}_{\s_0}(x_0,r)=(a,b).$}}
\label{simple metrics}
\end{figure}

\begin{definition}
\label{d1.1}
Let $X$ be an arbitrary nonempty set. \bd{Metric}\footnote{\text{Symbolize $d,\rho$ or $\sigma$.}} (or distance) in $X,$ is a map $\rho:X\times X\longrightarrow \R$ obeying the following properties:
\begin{enumerate}
\item  $\rho(x,y) \geq0,\ \forall x,y \in X$ and $\rho(x,y)=0 \Leftrightarrow x=y$
\item  $\rho(x,y)=\rho(y,x),\   \forall x,y \in X,$ (Symmetric property)
\item  $\rho(x,y)\leq \rho(x,z)+\rho(z,y),\   \forall x,y,z \in X,$ (Triangular inequality)
\end{enumerate}
The elements of the set are called \bd{points}, the real nonnegative number $\rho(x,y)$ is called the \bd{distance} between $x,y \in X$ and  the pair $(X,\rho)$ \bd{metric space}.
\end{definition}
Consequently, a set equipped with a metric, automatically obtains the structure of a topological space \footnote{\text{A topological space doesn't have to be a metric space.}}. We now define the open and the closed ball of center $x_0\in X$ and radius $r>0,$ notions necessary for the topological description of a metric space. 
\begin{definition}
\label{d1.2}
Let $(X,d)$ be a metric space, $x_0\in X$ and $r>0$. The set
$S_{d}(x_0, r)=\{x\in X: d(x, x_0)<r \}$ is called an \bd{open ball} of center $x_0$ and radius $r$.
\end{definition}
\begin{definition}
\label{d1.3}
Let $(X,d)$ be a metric space, $x_0\in X$ and $r\geq 0$. The set
$ \widetilde{S}_{d}(x_0, r)=\{x\in X: d(x, x_0)\leq r \}$ is called a \bd{closed ball} of center $x_0$ and radius $r$.
\end{definition}
\begin{definition}
\label{d1.4}
The set $A\subseteq X$ is called an \bd{open set}, if for every $\alpha \in A$ there exists $r>0$ such that $S_{d}(\alpha ,r)\subseteq A$.
\end{definition}
\begin{definition}
\label{d1.5}
The set $B\subseteq X$ is called a \bd{closed set}, if its complement $X\setminus B$ is an open set.
\end{definition}
\begin{definition}
\label{d1.6}
The set $\Gamma$ $\subseteq X$ is \bd{bounded} if there exists $x_{0} \in X$ and $r>0,$ such that $\Gamma \subseteq S_{d}(x_{0} ,r).$  
\end{definition}
\paragraph{Examples of Metric spaces:}
\begin{itemize}
\item The most common metrics to use in $\R^n$ are $d_1,\ d_2$ between its points $\mbd{x}=(x_1,x_2,...,x_n),$  $\mbd{y}=(y_1,y_2,...,y_n).$\\

The \bd{metric} $\mbd{d_{1}}$ (Manhattan metric) in $\R^n$ is defined as \[ d_{1}(\mbd{x}, \mbd{y})= \sum_{i=1}^n |x_{i}-y_{i}|.\]
Thus, in $\R^{2}$:  $ d_{1}(\mbd{x}, \mbd{y})= |x_{1}-y_{1}|+|x_{2}-y_{2}|$ and the closed ball is respectively $\widetilde{S}_{d_1}(\mbd{0}, r)=\{ \mbd{x}\in \R^2: d_{1}(\mbd{x, 0})\leq r \}$,  \big(Fig.~\ref{simple metrics} -  i), with $r=1$ \big).\\

The \bd{metric} $\mbd{d_{2}}$ (Euclidean metric) in $\R^n$ is defined as
 \[d_{2}(\mbd{x}, \mbd{y})= \Big( \sum_{i=1}^n(x_{i}-y_{i})^{2} \Big) ^{ \frac{1}{2} }. \] Thus, in $\R^2$: $d_2(\mbd{x,y})= \big( (x_1-y_1)^2+(x_2-y_2)^2 \big)^{\frac{1}{2}}$ and the closed ball is respectively $\widetilde{S}_{d_2}(\mbd{0}, r)=\{ \mbd{x}\in \R^2: d_{2}(\mbd{x, 0})\leq r \},$  \big(Fig.~\ref{simple metrics} -  ii), with $r=1$ \big).
\item In every nonempty set $X$ the \bd{metric} $\mbd{\s_{0}}$ (discrete metric) between points $x,y \in X$ is defined as
\begin{equation}
\label{1.1}
\s_{0}(x, y)=
\begin{cases}
 1, & \text { for $x\neq y,$}\\
 0, & \text { for $x = y$.}\\
\end{cases}
\end{equation}
Thus, the closed ball is $\widetilde{S}_{\s_0}(x_0, r)= \begin{cases} \{ x_0 \}, & 0<r<1\\ X, & r\geq 1 \end{cases},$  \big( Fig.~\ref{simple metrics} - iii), with $X=(\alpha,\beta)\subseteq \R$ \big).
\end{itemize}

\subsection{Cartesian product space}
If the set $X$ is arbitrary and not of a specific structure (e.g. vector space), the discrete metric (\ref{1.1}) seems to be the only available choice.

Given the metric spaces $(X_{i},\rho_{i}),\ i=1,2...,n$ we define the $p$-metrics ($p\geq 1$)\footnote{\text{For $0<p<1$ the triangular inequality does not hold, hence we do not define a metric.}} in the Cartesian product $X=X_{1}\times X_{2}\times...\times X_{n}$ for $\mathbf{x}=(x_1,...,x_n),\ \mathbf{y}=(y_1,...,y_n) \in X$:  
\begin{equation}
\label{1.2}
d_{(\rho_{1},\rho_{2},...,\rho_{n};\ p)} \big( \mathbf{x},\mathbf{y} \big)=
\begin{cases}
\Big( \sum_{i=1}^n \big( \rho_{i}(x_{i},y_{i})\big) ^{p} \Big) ^{\frac{1}{p}},\ & \text{for $1 \leq p<+ \infty $}\\
\max \{\rho_{i}(x_{i},y_{i}),\ i=1,2,...,n \}, & \text{for $p=+ \infty$}\\
\end{cases}
\end{equation} 
In case $X_{1}=X_{2}=...=X_{n}=Y$, i.e. $X=Y^n$ and $\rho_{1}=\rho_{2}=...=\rho_{n}= \rho$, we denote $d_{(\rho;\ p)}(\cdot\ ,\ \cdot)$ instead of $d_{(\rho_{1},\rho_{2},...,\rho_{n};\ p)}(\cdot\ ,\ \cdot).$

Metrics in (\ref{1.2}) are compatible with the ones that already exist in $X_{1},X_{2},...,X_{n}$ according to the following sense. Let $( y_{1},...,y_{n}) \in X$ be an arbitrary fixed point. Coinciding $x_{i}\in X_{i}$ with $( y_{1},...,x_{i},...,y_{n})\in X$ we have
\[d_{(\rho_{1},\rho_{2},...,\rho_{n};\ p)} \big( (y_{1},...,x_{i},...,y_{n}) ,\ \big( y_{1},...,x'_{i},...,y_{n}) \big)=\rho_{i}(x_{i},x'_{i}),\ \text{for}\  p\geq 1, \] where $i$ is the index corresponding to the metric space $(X_{i},\rho_{i}).$

\subsection{$\R^n\ (n\geq 1)$ equipped with metrics}
Due to the discrete nature of computers, our main interest is the set $X$ of the metric space to be a vectored space or subspace. In most of the applications the space that appears is $\R^{n}$ or subsets of this space. The axiomatic foundation of the set $\R $ of real numbers, gives us the latitude to define the metric $\s_{|\cdot |}(x,y)=|x-y|,\ x,\ y \in \R ,$ where $| \cdot |$ stands for the absolute value of a real number. More generally we may take the metrics $\s_{s}(x,y)=|x-y|^{s},$ for $0<s\leq 1\ (\s_{1}=\s_{| \cdot |}).$ At this point it is important to consider that
\begin{equation}
\label{1.3}
\lim_{s\rightarrow 0^{+}}\s_{s}(x,y)=\s_{0}(x,y).
\end{equation}

In case of the set $X\subseteq \R \times \R \times... \times \R=\R ^n$, emerge the $p$-metrics resulting from $(\R, \s_{| \cdot |})$, $(\R, \s_{s})$ for $0<s<1$ and $(\R, \s_{0})$ for points $\mathbf{x}=(x_1,...,x_n),\ \mathbf{y}=(y_1,...,y_n) \in \R^n,$ respectively. Alternatively, a combination of $\s$-metrics is also possible in order to measure in a different way among subsets of $X,$ however the latter choice lacks in practice.\\
Analytically we use the following metrics:
\paragraph{Usual metrics in $\R^n$:}
\begin{itemize}
\item Let $\R$ equipped with the metric $\s_{| \cdot |}=|x-y|.$ It follows that $\R^n$ is equipped with the metric $d_{(\s_{| \cdot |};\ p)}$ which according to (\ref{1.2}) leads to:
\begin{equation}
\label{1.4}
d_{(\s_{|\cdot |};\ p)} \big( \mathbf{x},\mathbf{y} \big)=
\begin{cases}
\Big( \sum_{i=1}^{n}|x_{i}-y_{i}|^{p} \Big) ^{\frac{1}{p}},\ & \text{for $1 \leq p<+ \infty $} \\
\max \{|x_{i}-y_{i}|,\ i=1,...,n \}, & \text{for $p=+ \infty $}\\
\end{cases}
\end{equation}
Thus, for $p=2$ we have the Euclidean metric in $\R^n$ (Fig.~\ref{pvariations}).

\begin{figure}[htb]
\centering%
\includegraphics[scale=0.6]{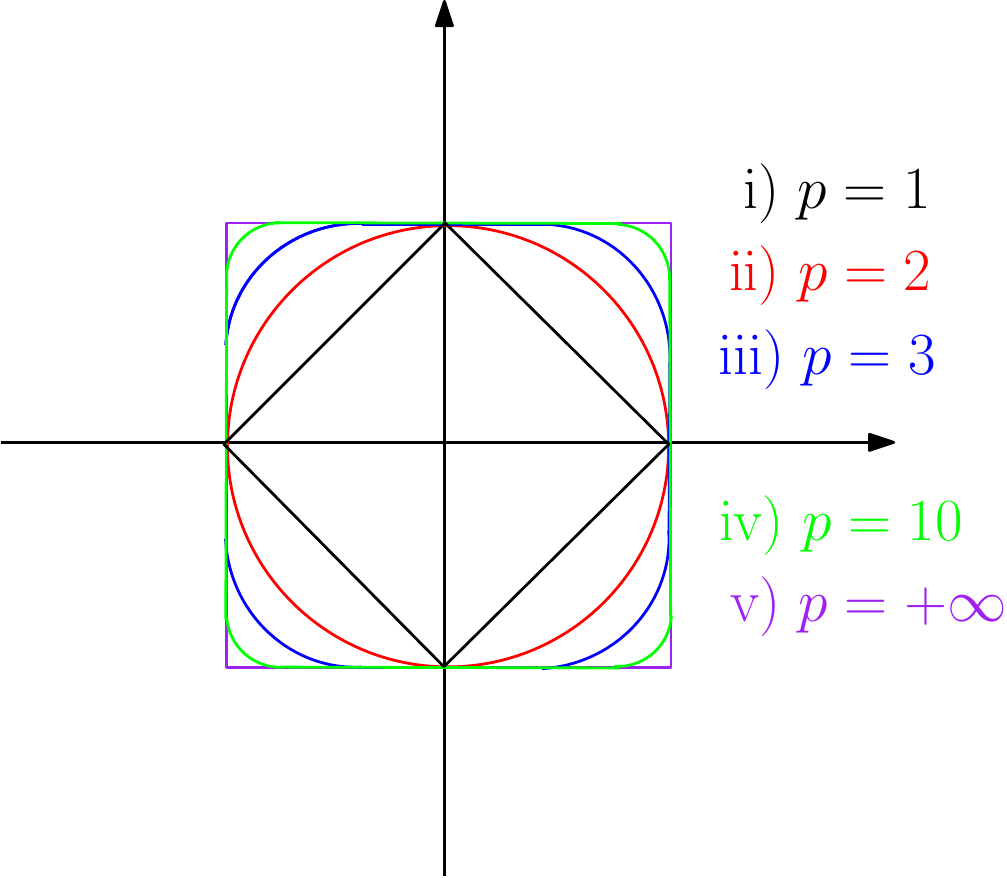}
\caption {\small{$\widetilde{S}_{d_{(\s_{| \cdot |};\ p)}}(\mbd{0},1)= \{ \mbd{x} \in \R^2 :\ d_{(\s_{| \cdot |};\ p)}(\mathbf{x},\mbd{0})\leq 1 \} $ for different values of $p$. Notice that while $p$ increases, the ball of our space tends to be $\widetilde{S}_{d_{(\s_{| \cdot |};\ +\infty )}}(\mbd{0},1),$ whereas $p$ decreases to $1,$ tends to be $\widetilde{S}_{d_{(\s_{| \cdot |};\ 1)}}(\mbd{0},1).$}}
\label{pvariations}
\end{figure}

\item Let $\R$ equipped with the metric $\s_{s}=|x-y|^{s}$ for $0<s\leq 1$. It follows that $\R^n$ is equipped with the metric $d_{(\s_{s};\ p)}$ that according to (\ref{1.2}) leads to:
\begin{equation}
\label{1.5}
d_{(\s_{s};\ p)} \big( \mathbf{x},\mathbf{y} \big)=
\begin{cases}
\Big( \sum_{i=1}^{n} \big( \s_{s}(x_{i},y_{i})\big) ^{p} \Big) ^{\frac{1}{p}},\ & \text{for $1 \leq p<+ \infty $} \\
\max \{\s_{s}(x_{i},y_{i}),\ i=1,...,n \}, & \text{for $p=+ \infty $}\\
\end{cases}
\end{equation} 
Specifically, for $p=1$ we have:
\begin{equation}
\label{1.6}
d_{(\s_{s};\ 1)} \big( \mathbf{x} ,\ \mathbf{y} \big)=\sum_{i=1}^n |x_i-y_i |^{s},\ 0<s\leq 1.
\end{equation} 
\end{itemize}
\paragraph{Discrete metric in $\R^n$:}\hspace{\fill}
\\
Let $\R$ equipped with the metric $\s_{0}(x, y)=
\begin{cases}
 1, & \text { for $x\neq y$}\\
 0, & \text { for $x = y$}\\
\end{cases}$, thus $\R^n$ is equipped with the metric $d_{(\s_{0};\ p)}$
that according to (\ref{1.2}) results to: 
\begin{equation}
\label{1.7}
d_{(\s_{0};\ p)} \big( \mathbf{x},\mathbf{y} \big)=
\begin{cases}
\Big( \sum_{i=1}^{n} \big( \s_{0}(x_{i},y_{i})\big) ^{p} \Big) ^{\frac{1}{p}},\ & \text{for $1 \leq p<+ \infty $} \\
\max \{\s_{0}(x_{i},y_{i}), \ i=1,...,n \}, & \text{for $p=+ \infty $}\\
\end{cases}
\end{equation}
Hence considering the case $p=1$ we have
\begin{equation}
\label{1.8}
d_{(\s_{0};\ 1)}(\mathbf{x} ,\ \mathbf{y})= \# \{ i: x_i \neq y_i \}.
\end{equation} 
As $d_{(\s_{s};\ 1)}$ is of most importance in sparse theory, we denote it as $d_s$, if  not to be confused with any other metric and use the symbolism $\widetilde{S}_{d_s}$ for the closed ball respectively. Finally, combining equation (\ref{1.8}), with both (\ref{1.3}) and (\ref{1.6}) we obtain:
\begin{equation}
\label{1.9}
\lim_{s \rightarrow 0^{+}}d_{s}(\mathbf{x},\mathbf{y})= \lim_{s \rightarrow 0^{+}} \Big ( \sum_{i=1}^n| x_i-y_i |^{s} \Big) =\# \{ i: x_i \neq y_i \}=d_{0}(\mathbf{x},\mathbf{y}).
\end{equation}

\begin{figure}[!]
\centering%
\includegraphics[scale=0.725]{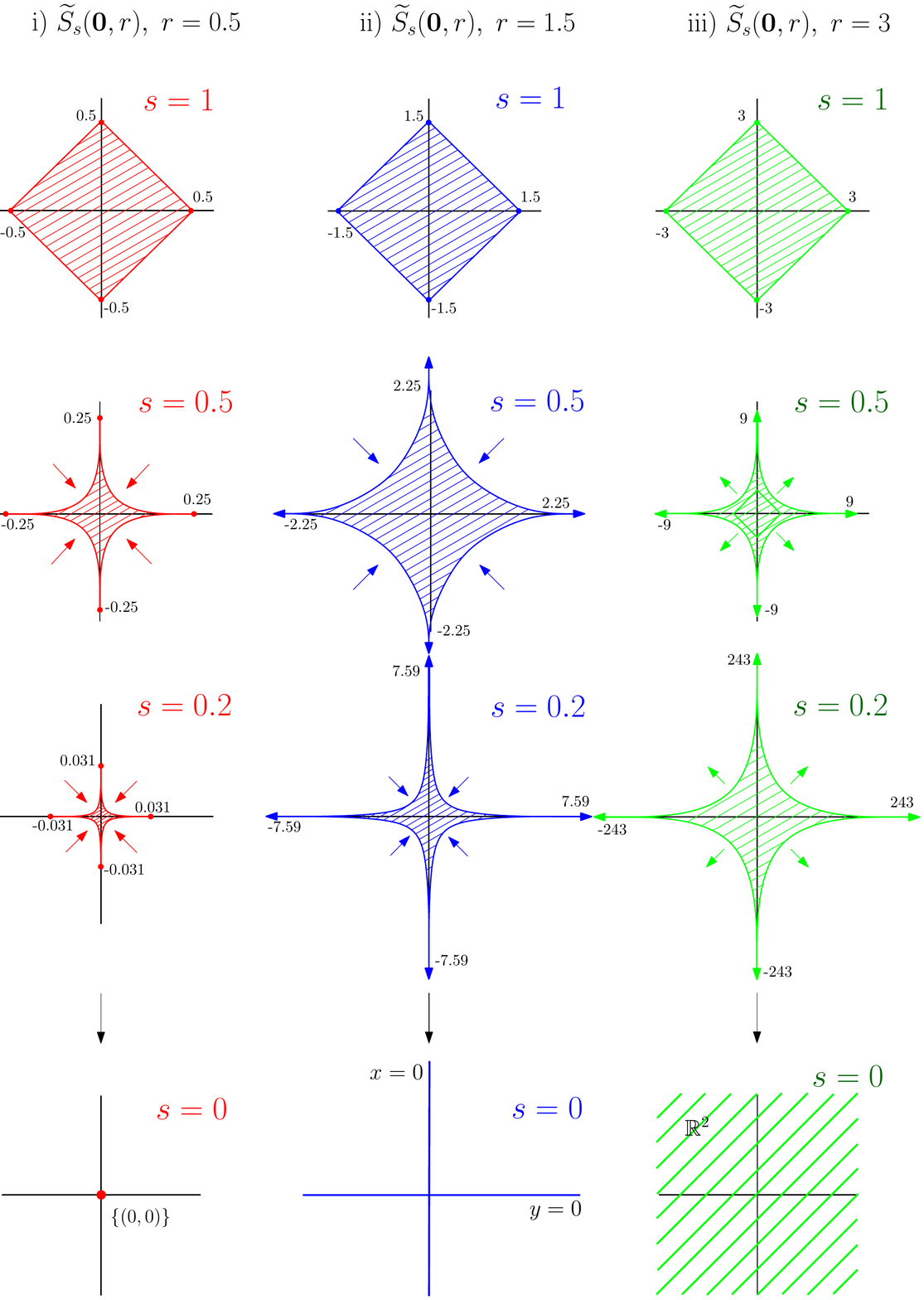}
\caption {\small{ $ \widetilde{S}_{s}(\mbd{0},r)= \{ \mbd{x}=(x,y) \in \R^2: \s_{s}(\mbd{x},\mbd{0})\leq r \},\ 0\leq s\leq 1 :$
\text{\en i)} For $0<r<1,\ \widetilde{S}_{0}(\mbd{0},r)=\{ (0,0) \}.$ 
\text{\en ii)} For $1 \leq r<2,\ \widetilde{S}_{0}(\mbd{0},r)=\{\mbd{x}=(x,y) \in \R^2: x=0$ or $y=0\}.$ 
\text{\en iii)} For $r\geq 2,\ \widetilde{S}_{0}(\mbd{0},r)=\R^2.$}}
\label{limits}
\end{figure}

The final equation indicates the behaviour of closed balls. In (Fig.~\ref{limits}) it can be easily seen that in $\R^2$ and for $r\in (0,1)$ ($r=0.5$) the balls $\widetilde{S}_{s}(\mbd{0},r)$ decrease, i.e. for $0<s'<s\leq 1$ we have $\widetilde{S}_{s'}(\mbd{0},r)\subset \widetilde{S}_{s}(\mbd{0},r)$ and finally tend to be the ball $\widetilde{S}_{0}(\mbd{0},r)=\bigcap \limits_{0<s\leq 1} \widetilde{S}_{s}(\mbd{0},r).$ For $r\in [1,2)$ ($r=1.5$) a relation of subset does not exist between the balls $\widetilde{S}_{s'}(\mbd{0},r)$ and $\widetilde{S}_{s}(\mbd{0},r)$ for $s'<s,$ however $\widetilde{S}_{s}(\mbd{0},r)$ decrease and tends to  coincide with the axis while $s\rightarrow 0^{+},$ i.e. the ball $\widetilde{S}_{0}(\mbd{0},r).$  For $r\in [2,+\infty)$ ($r=3$) the balls increase and for $0<s'<s\leq 1$ we have $\widetilde{S}_{s'}(\mbd{0},r)\supset \widetilde{S}_{s}(\mbd{0},r)$ until they finally fill the whole space, while $\widetilde{S}_{0}(\mbd{0},r)=\bigcup \limits_{0<s\leq 1} \widetilde{S}_{s}(\mbd{0},r).$

Equation (\ref{1.9}) indicates a desirable measure of sparsity, defined as
\begin{equation}
d_0(\mbd{x},\mbd{0})=\# \{ i: x_i\neq 0\},
\end{equation} 
measures the number of nonzero coordinates\footnote{Also called support of a vector and denoted as $supp\{\mbd{x}\}.$} of a vector and belongs to the family of metrics. 

\begin{figure}[htb]
\centering%
\includegraphics[scale=0.6]{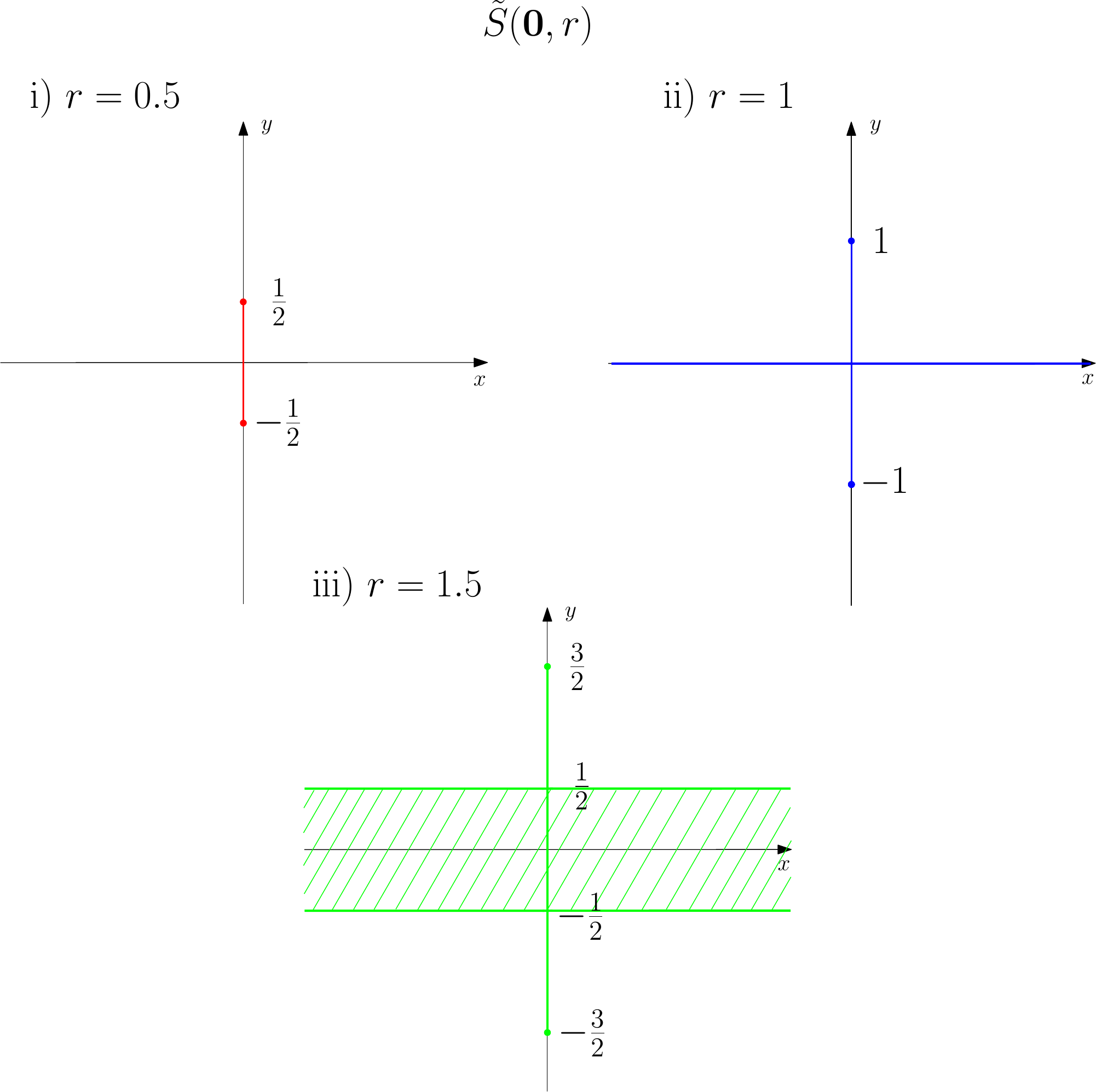}
\caption {\small{ $ \widetilde{S}(\mbd{0},r)=\{ \mbd{x}=(x,y)\in \R^2: \s_{0}(x,0)+ |y|\leq r \} $ i) For $0<r=0.5<1,\ \widetilde{S}(\mbd{0},0.5)=\{ -1/2\leq y\leq 1/2,\ x=0 \}.$ ii) For $r=1,\ \widetilde{S}(\mbd{0},1)=\{ -1\leq y\leq 1,\ x=0\ or\ y=0,\ x\neq 0\}.$ iii) For $r=1.5>1,\ \widetilde{S}(\mbd{0},1.5)=\{ -3/2 \leq y\leq 3/2,\ x=0\ or\ -1/2\leq y\leq 1/2,\ x\neq 0\}.$}} 
\label{altmetrics}
\end{figure}

\paragraph{Alternative metrics in $\R^n$:}\hspace{\fill}\\
Another measure constructed by a combination of different metrics (\ref{1.2}) enables us to measure each subset differently. At its simplest form we state an example in $\R^2.$ Let the metric space $ (\R,\s_{0}) \times (\R,\s_{|\cdot|})$ and set $\mbd{x}=(x_1,x_2),\ \mbd{y}=(y_1,y_2)\in \R^2.$ 
\[d_{(\s_0,\s_{|\cdot|};\ p)}(\mbd{x},\mbd{y})=\Big((\s_{0}(x_1,y_1))^{p}+ |x_2-y_2|^p \Big)^{1/p},\ 1\leq p< +\infty\]
Thus, for $p=1$:
\begin{equation}
\label{almetric1}
d_{(\s_0,\s_{|\cdot|};\ 1)}(\mbd{x},\mbd{y})=\s_{0}(x_1,y_1)+ |x_2-y_2|
\end{equation}
Consequently, the closed ball of center $\mbd{0}$ and radius $r$ are (Fig.\ref{altmetrics}):
\begin{equation}
\label{almetric2}
\widetilde{S}(\mbd{0},r)=\{ \mbd{x}\in \R^2: \s_{0}(x_1,0)+ |x_2|\leq r \}=\begin{cases} -r\leq x_2\leq r,\ & x_1=0,\\
1-r\leq x_2\leq r-1, & x_1\neq 0.
\end{cases}
\end{equation}

\hspace{\fill}
\\
\\
\\
\\
\\
\\
\\
\\
\\
\\
\\
\\
\\
\\

\section{Normed spaces}
 
\begin{definition}
\bd{Vector (linear) space} is called the trio ($X$, +, $\cdot$), where $X$ is a nonempty set,  $+\ :\ X \times X \longrightarrow X$ an inner operation (addition) and $\cdot\ : \ F\footnote{\text{The field $F= \mathbb{C}$ or $\R.$}}  \times X \longrightarrow X$ an outer operation (scalar product) that obey the following properties:
\begin{enumerate}
\item $x+y=y+x,\ \forall x,y \in X$,
\item $(x+y)+z=x+(y+z),\ \forall x,y,z \in X$,
\item There exists $0\in X$ such that $x+0= 0+x=x,\ \forall x\in X,$ 
\item For all $x\in X$ there exists $-x\in X$ such that $x+(-x)=(-x)+x=0,$
\item $\lambda (x+y)=\lambda x+ \lambda y,\ \forall x,y\in X$ and $\lambda \in F,$
\item $(\lambda + \mu)x=\lambda x+ \mu x,\ \forall x\in X$ and $\lambda,\ \mu \in F,$
\item $\lambda (\mu x)=(\lambda \mu)x,\ \forall x\in X$ and $\lambda,\ \mu \in F,$
\item $1 x=x,\ \forall x\in X.$
\end{enumerate}
The elements of a vectored space are called \bd{vectors}.
\label{d1.7}
\end{definition}

\begin{definition}
Let $X$ be a vector space over a field of numbers $F$. The set $A\subseteq X$ is called \bd{convex}, if for every pair $x,y\in X$ and every $t\in [0,1]$, the element $tx+(1-t)y$ belongs to the set $A$ as well.
\label{d1.8}
\end{definition}
\begin{definition}
A real function $f:A\rightarrow \R$ defined over a convex subset of a linear space $X$ is called \bd{convex}, if for every $x,y \in A$ and $t\in [0,1]$, 
\[f \big( tx+(1-t)y \big)\leq tf(x)+ (1-t)f(y). \] 
\label{d1.9}
\end{definition}
\begin{definition}
A real function  $f:A\rightarrow \R$ defined over a convex subset of a linear space $X$ is called \bd{concave}, if for every $x,y \in A$ and $t\in [0,1]$, 
\[f \big( tx+(1-t)y \big)\geq tf(x)+ (1-t)f(y). \] 
\label{d1.10}
\end{definition}

Let $\R$ be the vector space. Thus, the absolute value obeys the following properties:
\begin{enumerate}
\item $|x|\geq 0,$ $x\in \R$ and $|x|=0 \Leftrightarrow x=0$
\item $| \lambda x |=| \lambda | |x|,\ x\in \R,\ \lambda \in \R,$ (positive homogeneous)
\item $|x+y|\leq |x|+|y|,\ x,\ y \in \R,$ (triangular inequality)
\end{enumerate}
Therefore the function $f(x)= |x|$ is positive homogeneous, convex and $f(x)>0$ for $x\neq 0.$
A norm is the generalization of the absolute value in higher-dimensional vector spaces.
\begin{definition}
Let $(X, +, \cdot)$ be a real vector space. The map $\| \cdot \|:X\longrightarrow \R$ is called \bd{norm} if it obeys the following properties:
\begin{enumerate}
\item $\| x\| \geq0,\ \ \forall x\in X$ and $\| x\|=0\Leftrightarrow\ x=0,$ 
\item $\| \lambda x \|= | \lambda| \| x\|, \ \ \forall x\in X$ and $\lambda \in \R,$ (positive homogeneous)
\item $\|x+y\|\leq \|x\|+\|y\|,\ \ \forall x, y\in X,$ (triangular inequality)
\end{enumerate}
The pair $(X,\| \cdot \|$) is called a \bd{normed space}.
\label{d1.11}
\end{definition}
It follows that $f(x)=\|x\|$ is also a positive homogeneous, convex function with $f(x)>0$ for $x\neq 0.$ It is not difficult to see that if $\|x-y\|=f(x-y)=d(x,y)$ for $x,y\in X,$ then $d$ is a metric in $X$ with $d(x,0)=\|x\|.$ Moreover, if $\rho$ is a metric in a vector space $X$ satisfying the additional properties $\rho(x+z,y+z)=\rho(x,y),\ x,y\in X$ and $\rho(\lambda x,0)=|\lambda|\rho(x,0)$ with $x\in X,\ \lambda \in \R$ (positive homogeneous), then $\rho(x,0)=\|x\|$ is a norm. However, we will see that some of the metrics defined do not derive from norms.

Likewise metric spaces $X=X_{1}\times X_{2}\times...\times X_{n},$ the $p-$norms in $(X_{i},\| \cdot \|_{i})$ are defined.
\paragraph{Examples of normed spaces for $X=\R^n$:}
\begin{itemize}
\item The $\mbd{p-}$\bd{norms} for $1<p<+ \infty$: \[\|\mbd{x}\|_p= \Big( \sum_{i=1}^n|x_i|^p \Big)^{\frac{1}{p}}.\]
\item The \bd{Euclidean norm ($p=2$):}
\[\|\mbd{x}\|_{2}= \Big(\sum_{i=1}^n x_{i}^{2} \Big)^{\frac{1}{2}}.\]
\item The \bd{addition norm} ($\mbd{1-}$\bd{norm}) $\| \ \|_\mbd{1}$ and the \bd{norm} $\| \ \|_{\mbd{\infty}}$ respectively: \[\|\mbd{x}\|_1= \sum_{i=1}^n|x_i|,\] \[\|\mbd{x}\|_\infty =max\{|x_i|: i=1, 2,...,n \}. \]
\end{itemize}
At this point we should emphasise that $1\leq p\leq+ \infty,$ so that $\|\cdot \|_{p}$ is a norm in $\R^n$, which could be easily proved using the \text{\en Minkowski} inequality.

\begin{figure}[htb]
\centering%
\includegraphics[scale=0.33]{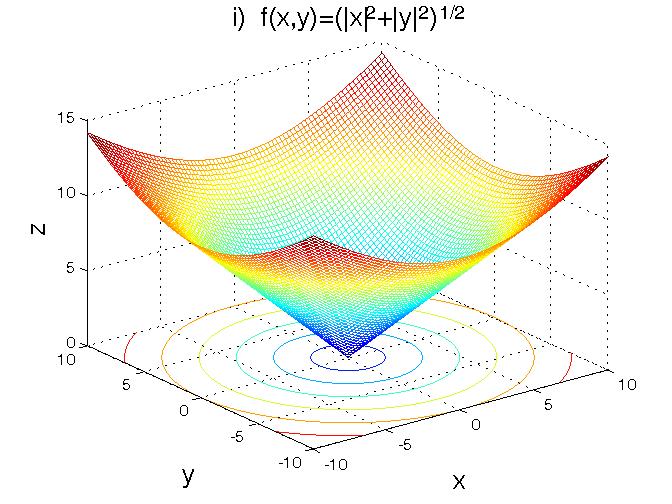}
\includegraphics[scale=0.33]{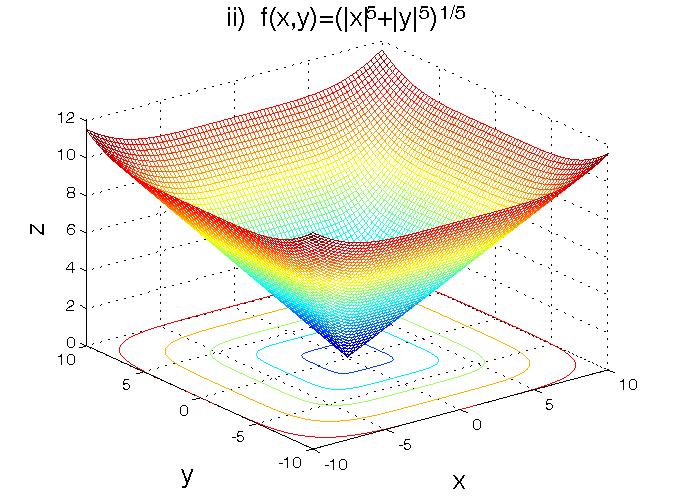}
\includegraphics[scale=0.33]{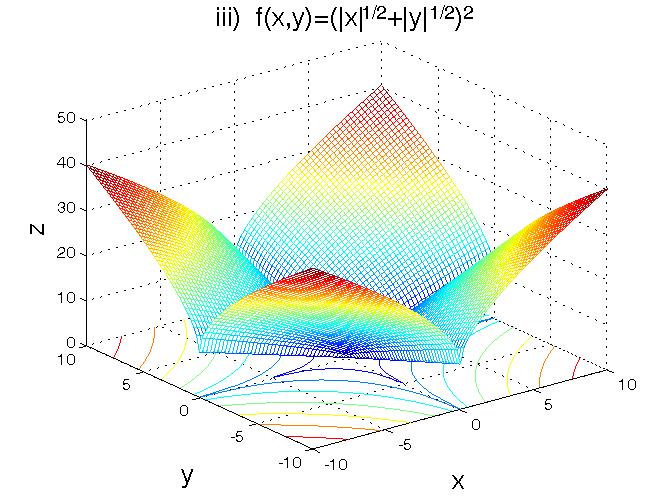}
\includegraphics[scale=0.33]{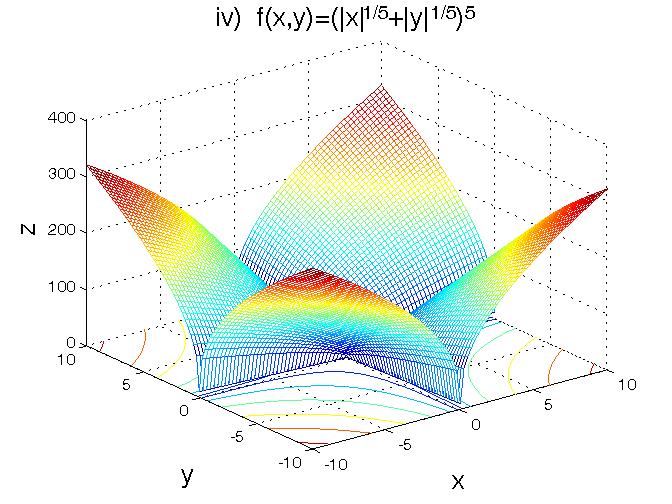}
\caption {\small{Function $f(x,y)=\big( |x|^p+|y|^p \big)^{\frac{1}{p}}$ for $p>0:$ \text{\en i) - ii)} For $p \geq 1$ those functions are convex. \text{\en iii) - iv)} For $0<p<1$ functions are not convex. \text{\en i) - iv)} All functions are positive homogeneous for all $p>0$.}}
\label{matlab1}
\end{figure}

Because of the demand for an optimization function, we set 
\begin{equation}
f_p(\mbd{x})=\Big( \sum_{i=1}^n|x_i|^p \Big)^{\frac{1}{p}},\ p>0
\label{1.10}
\end{equation}
which is positive homogeneous for every $p>0,$ whilst convex only for $p\geq 1.$ For $0<p<1$ the function is partially concave, hence the triangular inequality is not satisfied (Fig.~\ref{matlab1}).

\begin{figure}[!]
\centering%
\includegraphics[scale=0.33]{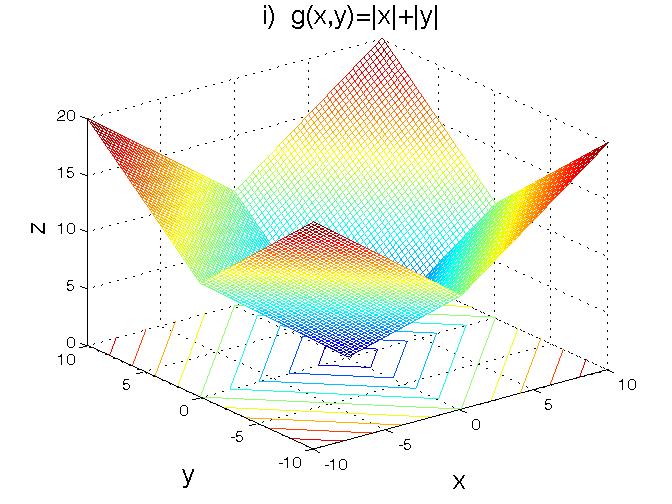}
\includegraphics[scale=0.33]{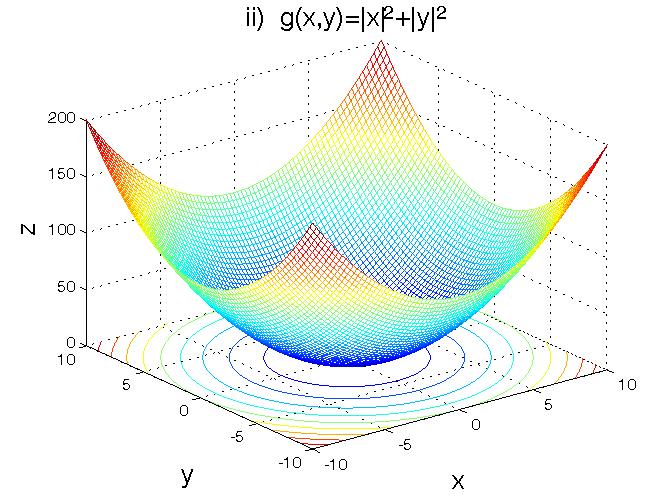}
\includegraphics[scale=0.33]{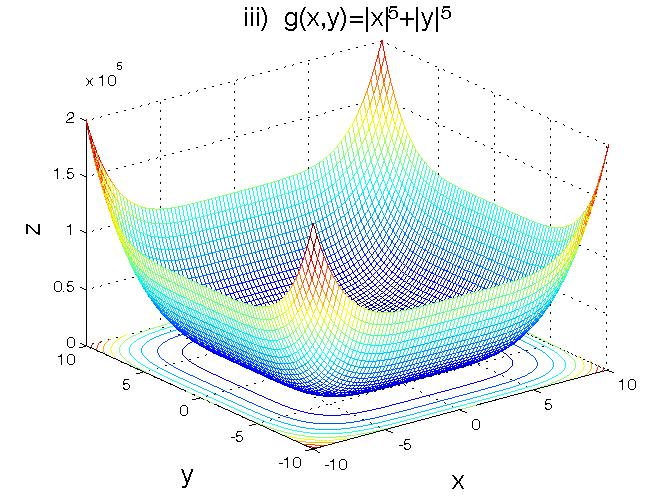}
\includegraphics[scale=0.33]{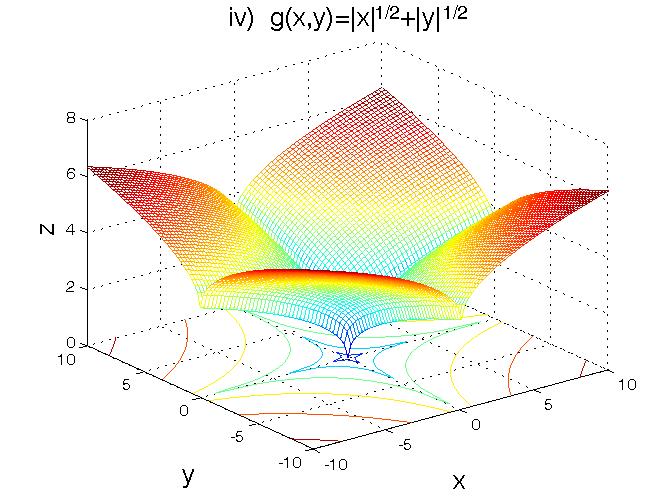}
\includegraphics[scale=0.33]{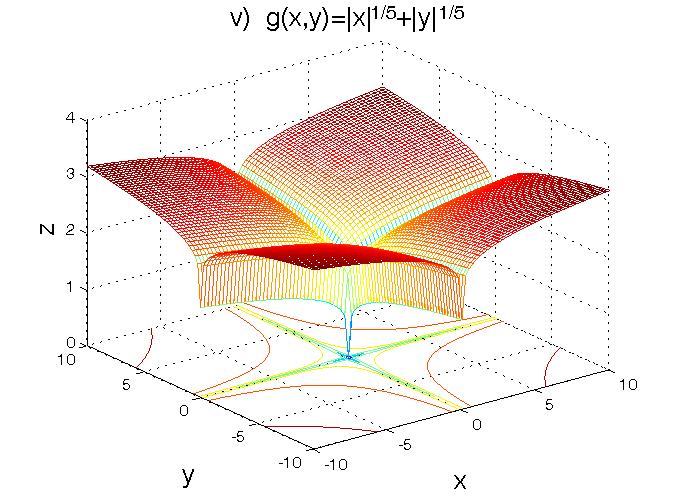}
\caption {\small{Function $g(x,y)=|x|^s+|y|^s$ for $s>0:$ \text{\en i) - ii)} For $s \geq 1$ the functions are convex. \text{\en iv) - v)} For $0<s<1$ functions are not convex. \text{\en ii) - v)} For $s\neq 1$ the functions are not positive homogeneous.}}
\label{matlab2}
\end{figure}

Metric $\s_{s}(x,y)=|x-y|^{s},\ x,\ y\in \R$ is transposition invariant, though not positive homogeneous for $s\neq 1.$ Hence, $d_{(\s_{s};\ p)}(\mbd{x},\mbd{0})$ for $0<s<1$ in $\R^n$ are not norms (Fig.~\ref{matlab2} for $p=1$).

From another point of view, the geometric interpretation gives us a clear image of all above. For $p\geq 1$ the closed balls are convex sets, unlike for $0<p<1$ .
Suppose $(X,\| \cdot \|)$ is a normed vector space, the balls $\widetilde{S}_{\|\cdot \|}(x,r)$ are always convex sets. Indeed, if $y,z \in \widetilde{S}_{\|\cdot \|}(x,r)$ then $\| y-x\|,\ \| z-x\|\leq r,$ thus for $\lambda \in (0,1),\ \|\lambda y+(1-\lambda z)-x\|=\|\lambda (y-x)+(1-\lambda )(z-x)\| \leq \lambda  r+(1- \lambda )r\leq r,$ i.e. $\lambda y +(1- \lambda)z \in \widetilde{S}_{\|\cdot \|}(x,r),$ hence the set $\widetilde{S}_{\|\cdot \|}(x,r)$ is convex.\\
\indt \bd{Remark}: The property of convexity is of great importance. Suppose that $(X,\| \cdot \|)$ is a normed vector space and let $K\subseteq X$ be a convex and symmetric ($K=-K$) open set, such that $R,r>0$ exist and $S_{ \| \cdot \| }(0,r) \subseteq K \subseteq S_{ \| \cdot  \|}(0,R).$ Thus $K$ constitutes the unit ball of another norm, i.e. $ \|x\|_{K}=inf\{ \lambda >0: x\in \lambda K \}$ (Minkowski functional) which is topological equivalent to the initial.\\
\\

\hspace{\fill}
\\
\\
\\
\\
\\
\\
\\
\\
\\
\\

\large{The University of Athens\\
Department of Mathematics\\
Panepistemiopolis 15784\\
Athens\\
Greece\\
Email:\\
ldalla@math.uoa.gr\\
ge99210@hotmail.com\\}

\end{document}